\documentclass[11pt,a4paper,fullpage]{article}

\usepackage{amsmath,amsfonts,amssymb,amsthm,graphics,amscd,mathtools,dsfont,float,hyperref,geometry}

\usepackage{fancybox}

\numberwithin{equation}{section}

\newtheorem{Teorema}{Theorem}[section]
\newtheorem{Propozicija}[Teorema]{Proposition}
\newtheorem{Definicija}[Teorema]{Definition}

\begin{document}

\title{On the fractional relaxation equation with Scarpi derivative}
\author{Matija Adam Horvat\footnote{University of Novi Sad, Faculty of Sciences, Novi Sad 21000, Serbia, e-mail: matija.adam.horvat@dmi.uns.ac.rs, ORCID ID 0009-0007-1850-4441}\footnote{Authors are supported by the Ministry of Education, Science and Technological Development of the Republic of Serbia under grants no. 451-03-33/2026-03/ 200125 and 451-03-34/2026-03/ 200125}\ \ and Nikola Sarajlija \ \footnote{corresponding author: Nikola Sarajlija, University of Novi Sad, Faculty of Sciences, Novi Sad 21000, Serbia; {\it e-mail}: {\tt nikola.sarajlija@dmi.uns.ac.rs}, ORCID id 0000-0003-2055-998X}}

\maketitle

\begin{abstract}
In this article we solve the Cauchy problem for the relaxation equation posed in a framework of variable order fractional calculus. After introducing some general mathematical theory we establish concepts of Scarpi derivative and transition functions which represent the essentials of our problem. Next, we provide an integral representation for the solution of our initial value problem where the transition function is chosen arbitrary. After that, we find an integral representation in several special cases in which we choose the transition function concretely. Finally, we give some numerical insights which prove our theoretical results.
\end{abstract}

\textit{$2020$ Math. Subj. Class:} 34A08, 34K37, 26A33, 34C26, 44A15.

\vspace{2mm}
\textit{Keywords and phrases:} fractional calculus, Scarpi derivative, relaxation equation.

\section{Introduction}

It is well known that an idea of non-integer differentiation and integration dates back to the work of Leibniz and Euler. However, the last 40 years have brought a real expansion in this area known as the fractional calculus. Namely, mathematicians realized that derivatives and integrals of non-integer order represent a very useful tool for modelling phenomena displaying non-localities in time. Therefore, many classical ordinary and partial differential equations have been successfully generalized to their fractional counterparts. For a detailed exposition of a study on constant-order differential and integral equations we suggest reference \cite{KNJIGA1}. \\

Soon after the interest in this topic encountered a huge boom, specialists noticed that in certain scenarios nature of non-localities can itself vary with respect to time. This is why the interest in constant order (CO) fractional operators partly moved to their variable order (VO) versions. During the years, there where many attempts at posing an appropriate theory for variable-order fractional calculus. We refer an interested reader to article \cite{PERSPEKTIVA} for more information. However, most of these approaches lacked in ri\-go\-rous mathematical framework and did not recover some basic properties of constant-order fractional operators. In this article we chose to rely on the variable-order theory proposed by an Italian engineer G. Scarpi in his paper \cite{SKARPI}, that experienced a recent revival in \cite{PERSPEKTIVA}. We found Scarpi's theory very appealing for our work, since the usual Caputo derivative turns out to be a special case of the Scarpi derivative and the fundamental theorem of calculus holds true for Scarpi derivative and integral. Thus, this kind of a variable-order generalization seemed the most natural for our purposes. Variable-order fractional calculus has been thoroughly investigated from a numerical perspective, among many papers see for example \cite{DODATNO1},\cite{DODATNO3},\cite{DODATNO2}, and there are only a few articles that do not rely on numerical machinery (\cite{DODATNO4}, \cite{DODATNO5}). In this work we invoke direct estimates rather than numerical ones, but for completeness we conduct a few numerical experiments at the end of the manuscript.\\

In this work we study the Cauchy problem for the VO fractional relaxation equation 
\begin{equation}\label{ScarpiRelaxation}
\begin{split}
\prescript{S}{}{D}_0^{\alpha(t)}u(t)=-\lambda u(t)\\
u(0)=u_0,
\end{split}
\end{equation}
where $\lambda>0$ is a relaxation parameter. In physics, value $-\lambda$ is called the elastic modulus. Symbol ${}^SD_0^{\alpha(t)}$ stands for Scarpi variable-order derivative and we introduce it in the next section. If ${}^SD_0^{\alpha(t)}$ were the usual derivative, we would have the analytic solution $u(t)=u_0e^{-\lambda t}$ of \eqref{ScarpiRelaxation}. If ${}^SD_0^{\alpha(t)}$ were the constant-order fractional derivative with order $\alpha\in(0,1)$, we would have the solution of the form $u(t)=u_0E_\alpha(-t^\alpha\lambda)$, where $E_\alpha(z)$ denotes the Mittag-Leffler function of order $\alpha$, that is $$E_\alpha(z)=\sum\limits_{k=0}^\infty\frac{z^k}{\Gamma(\alpha k+1)}.$$In this paper we solve the case that is the most general among mentioned, that is when ${}^SD_0^{\alpha(t)}$ is a variable order fractional derivative. Relaxation models are common throughout physics and natural sciences and are well studied, for example in papers \cite{VariableOrderRelaxation}, \cite{MainardiRelaxation}, \cite{Relaxation2}, \cite{Relaxation3}. We also mention an interesting article \cite{DistrOrderRelaxation} which presents a distributed order version of relaxation oscillations. Till this moment, equation (\ref{ScarpiRelaxation}) has been investigated only numerically in \cite{PERSPEKTIVA} and \cite{GarrappaNumerickoResenje}, where one can find some preliminary numerical computations. In this work, however, we provide an analytical formula for the solution to this problem. \\

The paper is organized as follows. In Section 2 we introduce the notation and some basic mathematical tools that we need in our work. Section 3 is our main section, where we find an explicit solution to problem \eqref{ScarpiRelaxation}. In Section 4 we illustrate our results on special types of transition functions. In the final section we provide an investigation of a result from Section 4 by numerical means.

\section{Mathematical prerequisites}

In the rest of this paper, we usually assume that functions in question have some decay properties related to integrability. Therefore, we will mostly exploit the concepts of Lebesgue integrable functions. We will need the concept of absolutely continuous functions (AC) on several occasions. Remember, $f\in AC[0,T]$ if $f$ is differentiable almost everywhere in $[0,T]$, $f'\in L^1[0,T]$ where $L^1[0,T]$ is the standard space of Lebesgue integrable functions, and $f(t)=f(0)+\int\limits_0^tf'(\tau)d\tau$ for all $ t\in[0,T].$ With $\mathcal{R}[f(t)]$ we denote the range space of function $f(t)$.\\

Since Scarpi's operators are introduced via characterizations in the Laplace domain, we now define the Laplace transform and state some of its important properties. Let $f$ be a piece-wise continuous function of exponential order $a$ defined on $[0,\infty]$, in other words inequality $\vert f(t)\vert\leq Me^{at}$ holds true for some positive $M$ and for all $t\geq t_0>0$. Then, for any $s\in\mathds{C}$ having $\mathit{Re}(s)>a$ the following indefinite integral converges:
$$F(s)=\mathcal{L}[f(t)](s)=\int\limits_0^\infty e^{-st}f(t) dt.$$
It is called the Laplace transform (LT) of function $f$. It is well known that such function $F(s)$ is unique, in a sense that if there is a piece-wise continuous function $f_1(t)$ of exponential order $b$ defined on $[0,\infty)$ and $\mathcal{L}[f]=\mathcal{L}[f_1]$, then $f(t)=f_1(t)$ on $[0,\infty)$ as an equality of piece-wise continuous functions. Hence, it is valid to define inverse Laplace transform, denoted by $\mathcal{L}^{-1}$, as
$$f(t)=\mathcal{L}^{-1}[F(s)](t)\ \Leftrightarrow\ F(s)=\mathcal{L}[f(t)](s). $$
The inverse Laplace transform is most commonly calculated using the Bromwich integral (also known as Mellin's inverse formula), that is
\begin{equation}\label{MELIN}
    \mathcal{L}^{-1} \left[ F(s) \right] (t) = \frac{1}{2 \pi i} \lim\limits_{T \to \infty} \int\limits_{\gamma - iT}^{\gamma + iT} F(s) e^{st} ds,
\end{equation}
where the integral is taken over line $Re(s) = \gamma$, where $\gamma$ is chosen such that all singularities of $F(s)$ have real part less then $\gamma$ and $F(s)$ is bounded on the line.

\noindent One can check that 
\begin{equation}\label{TAUBERIAN0}
\lim\limits_{\mathit{Re}(s)\to\infty}F(s)=0,
\end{equation}
and that the following Tauberian results are true.
\begin{Teorema}\cite{OPENHAJM}
For function $f$ as described above the following is true:
\begin{equation}\label{TAUBERIAN1}
\lim\limits_{t\to0+}f(t)=\lim\limits_{s\to\infty}s F(s),
\end{equation}
provided that $f'$ has an integrable singularity at $t=0$.\\
Suppose that every pole of $F(s)$ is either in the open left half plane or at the origin, and that $F(s)$ has at most a single pole at the origin. Then 
\begin{equation}\label{TAUBERIAN2}
\lim\limits_{t\to\infty}f(t)=\lim\limits_{s\to0}s F(s)
\end{equation}
\end{Teorema}

Let $\alpha\in(0,1)$. Standard Caputo derivative is defined by means of a convolution operator  \cite{KNJIGA1} as
\begin{equation}\label{KAPUTO}
\prescript{C}{}{D}_0^\alpha f(t)=\int\limits_0^t\varphi(t-\tau)f'(\tau)d\tau,\quad \varphi(t)=\frac{t^{-\alpha}}{\Gamma(1-\alpha)},  \quad t\geq0.  
\end{equation}
and it acts as the left inverse to the Riemann-Liouville fractional integral
\begin{equation}\label{RL}
\prescript{RL}{}{I}_0^\alpha f(t)=\int\limits_0^t\psi(t-\tau)f(\tau)d\tau,\quad \psi(t)=\frac{t^{\alpha-1}}{\Gamma(\alpha)},\quad t\geq0.
\end{equation}
It is easily deduced that 
\begin{equation}\label{LT}
\Phi(s)=\mathcal{L}[\varphi(t)](s)=s^{\alpha-1},\quad \Psi(s)=\mathcal{L}[\psi(t)](s)=s^{-\alpha}.
\end{equation}
\subsection{Transition function $\alpha(t)$}
In order to make a transition from constant order fractional operators to the ones having variable order, we have to replace constant order of differentiation $\alpha\in(0,1)$ in \eqref{KAPUTO} and \eqref{RL} by a certain function $\alpha(t)$ such that $\mathcal{R}[\alpha(t)]\subset(0,1)$. Some authors choose functions $\alpha$ having bounded domain, but we will not follow their steps in order to cover a more general case. Therefore, let $\alpha:[0,\infty)\to(0,1)$. This function $\alpha(t)$ will be considered as a variable order of differentiation, and thus we give it a special name. We call it the \emph{transition function}. It is easily concluded that function $\alpha:[0,\infty)\to(0,1)$ cannot be chosen arbitrarily. It must satisfy some necessary mathematical assumptions that will produce a suitable physical interpretation as well. Hence, we first assume that Laplace transformation of $\alpha(t)$ exists and its analytical expression is known. We denote it by $A(s)$, that is $A(s)=\mathcal{L}[\alpha(t)](s)$. \\

Next assumption on function $\alpha(t)$ that me must introduce is that there exists
\begin{equation}\label{PRETPOSTAVKA}
    \lim\limits_{t\to0+}\alpha(t)=\overline{\alpha}\in(0,1).
\end{equation}
Without this simple assumption, Laplace transform of the integral kernel appearing in definition of Scarpi derivative would not vanish at infinity, contrary to \eqref{TAUBERIAN0}. For more information on this phenomenon, consult \cite[p.7]{PERSPEKTIVA}. The assumption that the range of $\alpha(t)$ is within $(0,1)$ (sub-diffusive case) also agrees with the definition of Scarpi derivative as introduced in \cite{PERSPEKTIVA}. As is explained in Section 6 of this latter reference, it is possible to observe $\alpha(t)$ with range $(n,n+1)$ for some $n\in\mathds{N}$ (higher-order operators). We, however, will not pursue this point in this article. \\

Assumption \eqref{PRETPOSTAVKA} asserts that function $\alpha(t)$ admits some \emph{initial state value}. Potentially interesting physical scenarios demand that $\alpha(t)$ reaches some \emph{final state value} as well. Therefore, we next assume that function $\alpha(t)$ shows transition to some final state $\tilde{\alpha}$, in other words there exists
\begin{equation}\label{PRETPOSTAVKA2}
    \lim\limits_{t\to\infty}\alpha(t)=\tilde{\alpha}\in(0,1).
\end{equation}
This final state is sometimes only reached asymptotically as $t\to\infty$.\\

If transition is strictly monotone, like in the case of an exponential transition (\cite{GarrappaNumerickoResenje}), then constants in \eqref{PRETPOSTAVKA},\eqref{PRETPOSTAVKA2} obey $0<\overline{\alpha}<\tilde{\alpha}<1$ or $0<\tilde{\alpha}<\overline{\alpha}<1$ and order $\tilde{\alpha}$ is only reached asymptotically as $t\to\infty$.\\

To sum up, a certain function $\alpha:[0,\infty)\to(0,1)$ will be called a \emph{transition function} if and only if:
\begin{itemize}
    \item $A(s)=\mathcal{L}[\alpha(t)](s)$ exists and its analytical expression is known;
    \item $\lim\limits_{t\to0+}\alpha(t),\lim\limits_{t\to\infty}\alpha(t)$ exist and belong to $(0,1)$. 
\end{itemize}
\subsection{Definition of Scarpi's operators}
Notice that if $\alpha(t)$ is constant, put $\alpha(t)=\alpha$, then $A(s)=\frac{\alpha}{s}$ and relations (\ref{LT}) become
\begin{equation}\label{LT2}
\Phi_\alpha(s)=\mathcal{L}[\varphi(t)](s)=s^{sA(s)-1},\quad \Psi_\alpha(s)=\mathcal{L}[\psi(t)](s)=s^{-sA(s)}.
\end{equation}
Put $\varphi_\alpha(t)=\mathcal{L}^{-1}[\Phi_\alpha(s)](t)$ and $\psi_\alpha(t)=\mathcal{L}^{-1}[\Psi_\alpha(s)](t)$, where $\mathcal{L}^{-1}$ denotes inverse Laplace transform. Inspired by relations (\ref{LT2}) we arrive at definition of Scarpi's variable order fractional operators.
\begin{Definicija}\label{DEFINICIJA}\cite{PERSPEKTIVA}
Assume that $f$ is absolutely continuous. Then
\begin{equation}
    \prescript{S}{}{D}_0^{\alpha(t)}f(t)=\int\limits_0^t\varphi_\alpha(t-\tau)f'(\tau)d\tau,\quad t\in[0,\infty),
\end{equation}
is called Scarpi fractional derivative of variable order $\alpha(t)$.\\
Assume that $f$ is an integrable function. Then
\begin{equation}
    \prescript{S}{}{I}_0^{\alpha(t)}f(t)=\int\limits_0^t\psi_\alpha(t-\tau)f(\tau)d\tau,\quad t\in[0,\infty),
\end{equation}
is called Scarpi fractional integral of variable order $\alpha(t)$.
\end{Definicija}
When $\alpha(t)$ is constant ($\alpha(t)=\alpha\in(0,1)$), relations \eqref{LT2} reduce to \eqref{LT} and thus definitions of Scarpi's derivative and integral reduce to the usual definitions of Caputo derivative and Riemann-Liouville integral, respectively. Furthermore, Scarpi's operators introduced in Definition \ref{DEFINICIJA} satisfy the fundamental theorem of calculus (\cite{PERSPEKTIVA}):
\begin{equation}\label{FUNDAMENTALNATEOREMASKARPI}
    \prescript{S}{}{D}_0^{\alpha(t)}[\prescript{S}{}{I}_0^{\alpha(t)}f(t)]=f(t),\quad \prescript{S}{}{I}_0^{\alpha(t)}[\prescript{S}{}{D}_0^{\alpha(t)}f(t)]=f(t)-f(0).
\end{equation}
In the sequel we will need the following formula for Laplace transform of $\prescript{S}{}{D}_0^{\alpha(t)}f(t)$, that can be easily proved using some basic properties of LT:
\begin{equation}\label{SKARPIJEVLAPLAS}
\mathcal{L}[\prescript{S}{}{D}_0^{\alpha(t)}f(t)](s)=s^{sA(s)}F(s)-s^{sA(s)-1}f(0).    
\end{equation}
We want to emphasize that Scarpi's differential operators are superior to other variable-order fractional derivatives in various senses. For example, if we try to define variable-order differential operators by naive replacement of $\alpha\in(0,1)$ by $\alpha(t):[0,\infty)\to(0,1)$ in \eqref{KAPUTO} and \eqref{RL}, than $D_0^{\alpha(t)}$ obtained in this way does not necessarily act as the left inverse to $I_0^{\alpha(t)}$ (see \cite{SAMKO}). We could also try such a naive replacement to other well known fractional derivatives such as the Weyl or Hadamard Fractional Derivative (see \cite{POSLEDNJAKNJIGA}), but Weyl's derivative demands posing non-local and impractical initial conditions while Hadamard's derivative is poorly suited to initial-value problems at $t=0$. Another attempt of generalizing relations \eqref{KAPUTO} and \eqref{RL} is a variable-order Caputo–Fabrizio derivative as introduced in \cite{ZHENGWANGFU}, but its derivative of a constant function is always zero, regardless of order, which is physically restrictive. Finally, there are many definitions that appear very useful for physical applications but lack mathematical precision, such as those in \cite{COIMBRA} or \cite{LORENZO}. However, they face conceptual mathematical problems as described in \cite{SAMKO} or \cite{ROSS}. To summarize, Scarpi's operators are both mathematically precise and physically appealing, with their definition being rather simple, and considering that they recover many known properties from constant-order fractional calculus (such as the fundamental theorem of calculus or semi-group property \cite{PERSPEKTIVA}), we conclude that Scarpi's differential operators are by far the best choice for our needs. Therefore, we procede with their use in the sequel.

\section{Integral representation of the unique solution to the Cauchy problem}

Our goal in this section is to find an explicit expression for the solution of problem (\ref{ScarpiRelaxation}). We will prove the following theorem.
\begin{Teorema}
Let $u_0$ in (\ref{ScarpiRelaxation}) be given. Then, the unique solution of (\ref{ScarpiRelaxation}) is given by
\begin{equation}\label{konacnoResenje}
\begin{aligned}
     &u(t) = \frac{u_{0}}{2 \pi i} \left( \mathfrak{R} + \int\limits_{0}^{\infty} \frac{\rho^{- \rho A \left( \rho e^{- i \pi} \right) - 1} e^{i \pi \left( \rho A \left( \rho e^{- i \pi} \right) + 1 \right)}}{\rho^{- \rho A \left( \rho e^{- i \pi} \right)} e^{i \pi \rho A \left( \rho e^{- i \pi} \right)} + \lambda} e^{- \rho t} d\rho \right. \\ &{} \hspace{13em}\left. - \int\limits_{0}^{\infty} \frac{\rho^{- \rho A \left( \rho e^{i \pi} \right) - 1} e^{- i \pi \left( \rho A \left( \rho e^{i \pi} \right) + 1 \right)}}{\rho^{- \rho A \left( \rho e^{i \pi} \right)} e^{- i \pi \rho A \left( \rho e^{i \pi} \right)} + \lambda} e^{- \rho t} d\rho \right),
\end{aligned}\end{equation}
where $\mathfrak{R}$ denotes the sum of all residues of $\mathcal{L}[u](s)$ which are located inside $\Gamma$ multiplied with $2\pi i$, and $\Gamma$ is a contour sketched below.
\end{Teorema}
\textbf{Proof. }Since we have an obvious linearity in the equation, we apply Laplace transform to both sides of \eqref{ScarpiRelaxation}. One obtains
$$s^{sA(s)}\tilde{u}(s)-s^{sA(s)-1}u_0+\lambda\tilde{u}(s)=0,$$
and therefore
$$\tilde{u}(s)=\frac{s^{sA(s)-1}}{s^{sA(s)}+\lambda}u_0,$$
where tilde  $\tilde{}$  denotes the LT of a function.\\

Next, we aim to compute the solution $u(t)=\mathcal{L}^{-1}[\tilde{u}(s)](t)$ using formula (\ref{MELIN}). It is clear that this formula gives us
\begin{equation}\label{RESENJE}
u(t)=\mathcal{L}^{-1} \left[ \tilde{u}(s) \right] (t) = u_0\cdot\frac{1}{2 \pi i} \lim\limits_{T \to \infty} \int\limits_{s_{0} - iT}^{s_{0} + iT} \frac{s^{sA(s)-1}}{s^{sA(s)}+\lambda}e^{st} ds,
\end{equation}
where $s_{0}>0$ is chosen such that all singularities of $\tilde{u}(s)$ have real part less then $s_{0}$. If we look at the contour given below, we see that inverse Laplace transform of $\tilde{u}(s)$ is identified as an integral over line $\Gamma_{v}$ when we let $R \to \infty$.
\begin{figure}[H]
    \centering
    \includegraphics[width=0.6\linewidth]{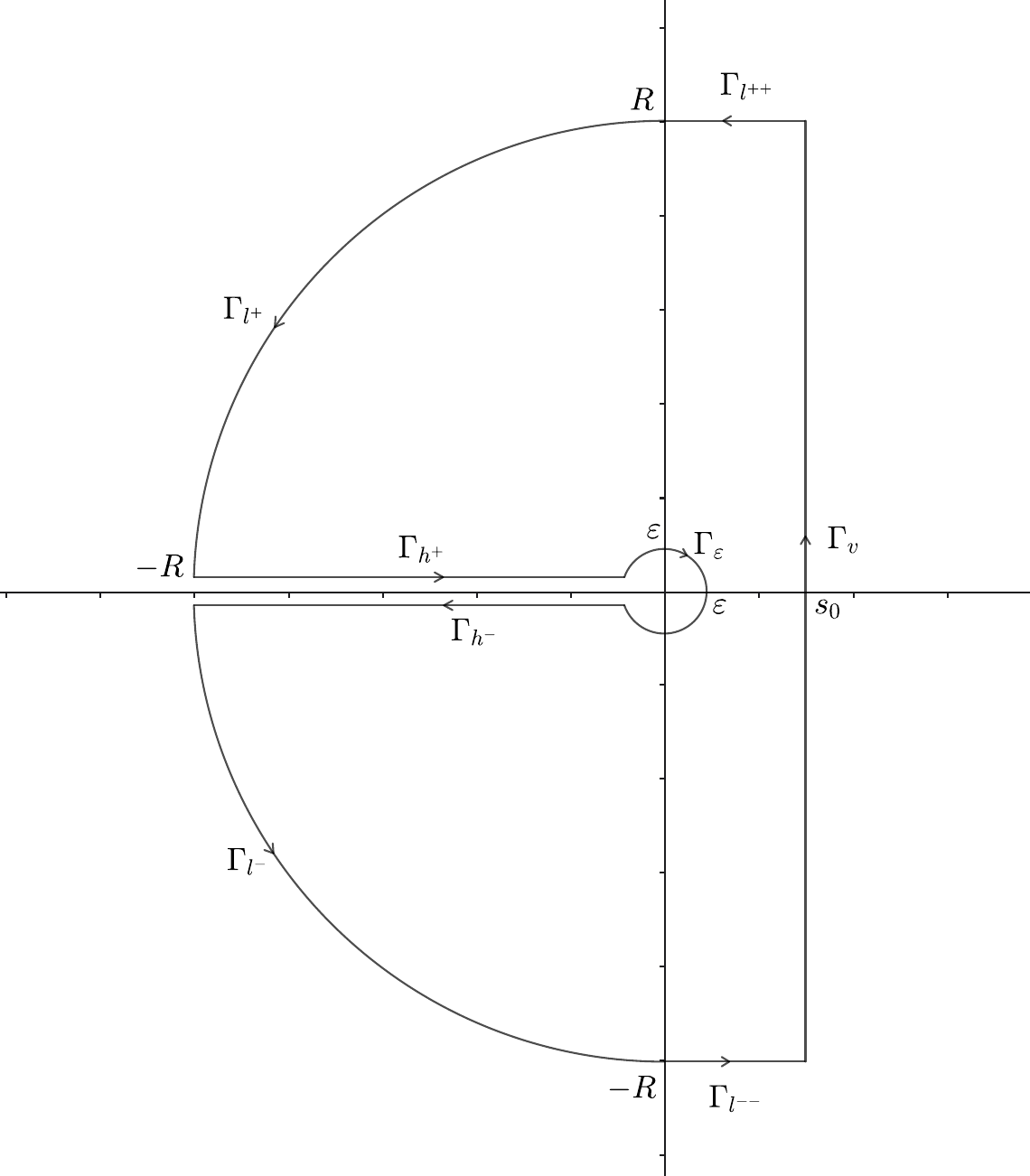}
    \caption{Contour of integration $\Gamma$}
    \label{kontura}
\end{figure}

We denote this contour by $\Gamma$, that is $$\Gamma=\Gamma_{l^{++}}\cup\Gamma_{l^+}\cup\Gamma_{h^+}\cup\Gamma_{\varepsilon}\cup\Gamma_{h^-}\cup\Gamma_{l^-}\cup\Gamma_{l^{--}}\cup\Gamma_v,$$
where, for arbitrarily chosen $R>0$ and $0<\varepsilon<R$ we have:
\begin{description}\label{OPIS}
    \item[$\Gamma_{l^{++}}$]: $s=-\xi s_0+iR,\ -1<\xi<0;$ 
    \item[$\Gamma_{l^{+}}$]: $s=Re^{i\varphi},\ \frac{\pi}{2}<\varphi<\pi$;
    \item[$\Gamma_{h^{+}}$]: $s=e^{i \pi}\xi,\ -R<-\xi<-\varepsilon$;
    \item[$\Gamma_{\varepsilon}$]: $s=\varepsilon e^{i\varphi},\ -\pi<\varphi<\pi$;
    \item[$\Gamma_{h^{-}}$]: $s=e^{- i \pi}\xi,\ \varepsilon<\xi<R$;
    \item[$\Gamma_{l^{-}}$]: $s=Re^{i\varphi},\ -\pi<\varphi<-\frac{\pi}{2}$;
    \item[$\Gamma_{l^{--}}$]: $s=\xi s_0-Ri,\ 0<\xi<1$;
    \item[$\Gamma_v$]: $s=s_{0}+\xi i,\ -R<\xi<R$. 
\end{description}

By the Cauchy residue theorem we have that 
\begin{equation*}
    \oint\limits_\Gamma\tilde{u}(s)e^{st}ds=\mathfrak{R},
\end{equation*}
where $\mathfrak{R}$ denotes the sum of all residues which are located inside $\Gamma$ multiplied with $2\pi i$. Next, we divide
\begin{equation*}
\oint_\Gamma\tilde{u}(s)e^{st}ds=\int_{\Gamma_{l^{++}}}+\int_{\Gamma_{l^{+}}}+\int_{\Gamma_{h^{+}}}+\int_{\Gamma_{\varepsilon}}+\int_{\Gamma_{h^{-}}}+\int_{\Gamma_{l^{-}}}+\int_{\Gamma_{l^{--}}}+\int_{\Gamma_v},
\end{equation*}
where it is understood that $\oint\limits_\Gamma\tilde{u}(s)e^{st}ds$ depends on $R$ and $\varepsilon$.
In what follows, we prove that when $R \to \infty$ and $\varepsilon \to 0+$ one gets the following equality: 
\begin{equation}\label{TRAZENO}\begin{aligned}
    \lim\limits_{\substack{R \to \infty \\ \varepsilon \to 0+}} \left( \int\limits_{\Gamma_{v}} \frac{s^{s A(s) - 1}}{s^{s A(s)} + \lambda} e^{st} ds + \int\limits_{\Gamma_{h^{+}}} \frac{s^{s A(s) - 1}}{s^{s A(s)} + \lambda} e^{st} ds + \int\limits_{\Gamma_{h^{-}}} \frac{s^{s A(s) - 1}}{s^{s A(s)} + \lambda} e^{st} ds \right) = \mathfrak{R}.
\end{aligned}\end{equation}
In other words, we prove that
\begin{equation*}\begin{aligned}
    \lim\limits_{\substack{R \to \infty \\ \varepsilon \to 0+}}\left( \int\limits_{\Gamma_{l^{++}}}\tilde{u}(s)e^{st}ds+\int\limits_{\Gamma_{l^{+}}}\tilde{u}(s)e^{st}ds + \int\limits_{\Gamma_{\varepsilon}}\tilde{u}(s)e^{st}ds  +\int\limits_{\Gamma_{l^{-}}}\tilde{u}(s)e^{st}ds+\int\limits_{\Gamma_{l^{--}}}\tilde{u}(s)e^{st}ds\right) = 0.
\end{aligned}\end{equation*}
%KOMENTAR!! Ovde bismo voleli da je $\mathfrak{R} = 0$ tj. da nekako uspemo da poka\v zemo da funkcija nema polove unutar konture, ovo mislim da \' ce nam najte\v ze biti. Kakogod i da nije nula imamo neku integralnu reprezentaciju inverzne Laplasove za to \v sto nam treba uz neki jo\v s \v clan ako ovo ne poka\v zemo !!

First, we prove that $\displaystyle \lim\limits_{R\to\infty}\int\limits_{\Gamma_{l^+}}\frac{s^{sA(s)-1}}{s^{sA(s)}+\lambda}e^{st}ds=0$. Having in mind definition of $\Gamma_{l_+}$ one gets
$$\int\limits_{\Gamma_{l^+}}\frac{s^{sA(s)-1}}{s^{sA(s)}+\lambda}e^{st}ds=\int\limits_{\frac{\pi}{2}}^\pi\frac{e^{(Re^{i\theta}A(Re^{i\theta})-1)Ln(Re^{i\theta})}}{e^{Re^{i\theta}A(Re^{i\theta})Ln(Re^{i\theta})}+\lambda}iRe^{i\theta}e^{Re^{i\theta}t}d\theta.$$
Evaluating the complex modulus of the latter integral and using the inverse triangle inequality in the denominator one obtains
$$\lim\limits_{R\to\infty} \left| \int\limits_{\Gamma_{l^+}}\frac{s^{sA(s)-1}}{s^{sA(s)}+\lambda}e^{st}ds \right| \leq\lim\limits_{R\to\infty}\int\limits_{\frac{\pi}{2}}^\pi \frac{\left|R^{Re^{i\theta}A(Re^{i\theta})}\right|}{\left|R^{Re^{i\theta}A(Re^{i\theta})}\right|-\lambda} e^{Rt\cos\theta} d\theta.$$
Since it is clear that $\lim\limits_{R\to\infty}Re^{i\theta}A(Re^{i\theta})=\lim\limits_{s\to\infty}sA(s)=\overline{\alpha}\in(0,1)$ by (\ref{TAUBERIAN1}) and that $t\cos\theta<0$ for $t>0,\ \theta\in(\frac{\pi}{2},\pi)$, we conclude that the latter limit is equal to 0.\\

Similar reasoning is deduced for the integral of $\tilde{u}(s)e^{st}$ over $\Gamma_{l^-}$, in other words
$$\lim\limits_{R\to\infty}\int\limits_{\Gamma_{l^-}}\frac{s^{sA(s)-1}}{s^{sA(s)}+\lambda}e^{st}ds=0$$
holds, and the proof is analogous to the one that was just seen.\\

Analogous result holds for the integrals of $\tilde{u}(s)e^{st}$ over segments $\Gamma_{l^{++}}$ and $\Gamma_{l^{--}}$. Namely, we shall exhibit that 
\begin{equation}\label{DVALIMESA}
\displaystyle \lim\limits_{R \to \infty} \int\limits_{\Gamma_{l^{++}}} \frac{s^{s A(s) - 1}}{s^{s A(s)} + \lambda} e^{st} ds = 0,\quad \displaystyle \lim\limits_{R \to \infty} \int\limits_{\Gamma_{l^{--}}} \frac{s^{s A(s) - 1}}{s^{s A(s)} + \lambda} e^{st} ds = 0.
\end{equation}
Observe that $\Gamma_{l^{++}}$ is exactly the line segment joining the points $s_{0} + iR$ and $iR$ (Fig. \ref{kontura}). Thus, by its definition, the following is true:
\begin{equation*}
    \int\limits_{\Gamma_{l^{++}}} \frac{s^{s A(s) - 1}}{s^{s A(s)} + \lambda} e^{st} ds = - \int\limits_{0}^{1} \frac{\left( \theta s_{0} + iR \right)^{\left( \theta s_{0} + iR \right) A\left( \theta s_{0} + iR \right) - 1} }{\left( \theta s_{0} + iR \right)^{\left( \theta s_{0} + iR \right) A\left( \theta s_{0} + iR \right)} + \lambda} e^{(\theta s_{0} + iR)t} s_0d\theta.
\end{equation*}
Taking the complex module in the equality above and evaluating the right hand side, simple calculations give us that
\begin{equation*}
    \left| \int\limits_{\Gamma_{l^{++}}} \frac{s^{s A(s) - 1}}{s^{s A(s)} + \lambda} e^{st} ds \right| \leq \int\limits_{0}^{1} \frac{s_{0} e^{\theta s_{0}}}{\sqrt{\theta^2s_0^2+R^2} \left| 1 + \frac{\lambda}{\left( \theta s_{0} + iR \right)^{\left( \theta s_{0} + iR \right) A\left( \theta s_{0} + iR \right)}} \right|} d\theta
\end{equation*}
Now, letting $R \to \infty$ and using \eqref{TAUBERIAN1} one gets that
\begin{equation}\label{Gama_l++}
    \lim\limits_{R \to \infty} \int\limits_{\Gamma_{l^{++}}} \frac{s^{s A(s) - 1}}{s^{s A(s)} + \lambda} e^{st} ds = 0.
\end{equation}
Thus, we have proved the first equality in \eqref{DVALIMESA}, and the second is proved similarly.

We now go on to show that the integral of $\tilde{u}(s)e^{st}$ over $\Gamma_{\varepsilon}$ tends to $0$ as $\varepsilon \to 0+$. Keeping in mind definition of $\Gamma_{\varepsilon}$, we have
\begin{equation*}
    \int\limits_{\Gamma_{\varepsilon}} \frac{s^{s A(s) - 1}}{s^{s A(s)} + \lambda} e^{s t} ds =  \int\limits_{- \pi}^{\pi} \frac{\left( \varepsilon e^{i \theta} \right)^{\varepsilon e^{i \theta} A \left( \varepsilon e^{i \theta} \right) - 1}}{\left( \varepsilon e^{i \theta} \right)^{\varepsilon e^{i \theta} A \left( \varepsilon e^{i \theta} \right)} + \lambda} e^{\varepsilon e^{i \theta} t}\cdot i \varepsilon e^{i \theta} d\theta.
\end{equation*}
Since it is clear that $\lim\limits_{\varepsilon\to0+}\varepsilon e^{i\theta}A(\varepsilon e^{i\theta})=\lim\limits_{s\to0}sA(s)=\tilde{\alpha}\in(0,1)$ by (\ref{TAUBERIAN2}), we conclude that the latter integral tends of to $0$ as we let $\varepsilon \to 0+$.
Hence, we have shown that relation (\ref{TRAZENO}) holds.\\

Having in mind definition of segments $\Gamma_{h^{+}}$ and $\Gamma_{h^{-}}$ and the fact that $s^{s A(s)} = e^{s A(s) Ln(s)}$, $s^{s A(s) - 1} = e^{(s A(s) - 1) Ln(s)}$, we write out the integrals over these segments as $R \to \infty$ and $\varepsilon \to 0+$ keeping in mind what branch of the complex logarithm we are taking for each segment. We obtain:
\begin{equation}\label{horizontalniIntegrali}\begin{aligned}
    &\lim\limits_{\substack{R \to \infty \\ \varepsilon \to 0+}}\int\limits_{\Gamma_{h^{+}}} \frac{s^{s A(s) - 1}}{s^{s A(s)} + \lambda} e^{st} ds = \int\limits_{0}^{\infty} \frac{\rho^{- \rho A \left( \rho e^{i \pi} \right) - 1} e^{- i \pi \left( \rho A \left( \rho e^{i \pi} \right) + 1 \right)}}{\rho^{- \rho A \left( \rho e^{i \pi} \right)} e^{- i \pi \rho A \left( \rho e^{i \pi} \right)} + \lambda} e^{- \rho t} d\rho, \\
    &\lim\limits_{\substack{R \to \infty \\ \varepsilon \to 0+}}\int\limits_{\Gamma_{h^{-}}} \frac{s^{s A(s) - 1}}{s^{s A(s)} + \lambda} e^{st} ds = - \int\limits_{0}^{\infty} \frac{\rho^{- \rho A \left( \rho e^{- i \pi} \right) - 1} e^{i \pi \left( \rho A \left( \rho e^{- i \pi} \right) + 1 \right)}}{\rho^{- \rho A \left( \rho e^{- i \pi} \right)} e^{i \pi \rho A \left( \rho e^{- i \pi} \right)} + \lambda} e^{- \rho t} d\rho.
\end{aligned}\end{equation}
Next, we prove convergence of the integrals in \eqref{horizontalniIntegrali}. We will provide a proof for the integral over $\Gamma_{h^{+}}$, and by analogy one can carry out a proof for the integral over $\Gamma_{h^{-}}$. 
Using basic calculations, we first write the integral over $\Gamma_{h^+}$ in the following form:
\begin{equation*}
    \lim\limits_{\substack{R \to \infty \\ \varepsilon \to 0+}}\int\limits_{\Gamma_{h^{+}}} \frac{s^{s A(s) - 1}}{s^{s A(s)} + \lambda} e^{st} ds = \int\limits_{0}^{\infty} \frac{e^{- \rho t}}{\rho e^{i \pi} \left( 1 + \frac{\lambda}{\rho ^{\rho e^{i \pi} A \left( \rho e^{i \pi} \right)} \left( e^{i \pi} \right)^{\rho e^{i \pi} A \left( \rho e^{i \pi} \right)}} \right)} d\rho.
\end{equation*}
And now, we examine the asymptotic behavior of the function under the integral sign as $\rho \to \infty$. Namely, recalling \eqref{TAUBERIAN1} and \eqref{PRETPOSTAVKA} we arrive at $\rho e^{i \pi} A \left( \rho e^{i \pi} \right) \to \overline{\alpha}$ as $\rho \to \infty$. Thus, the asymptotic behavior of the function under the integral as $\rho \to \infty$ can be described as 
\begin{equation}\label{roUBeskonacno}
    \frac{e^{- \rho t}}{\rho e^{i \pi} \left( 1 + \frac{\lambda}{\rho ^{\rho e^{i \pi} A \left( \rho e^{i \pi} \right)} \left( e^{i \pi} \right)^{\rho e^{i \pi} A \left( \rho e^{i \pi} \right)}} \right)} \sim \frac{- e^{- \rho t}}{\rho + \lambda e^{- i \pi \overline{\alpha}} \rho^{1 - \overline{\alpha}}}, \quad \rho \to \infty.
\end{equation}
Next, we turn to the asymptotic behavior of the same function as $\rho \to 0$, and now invoking \eqref{TAUBERIAN2} and \eqref{PRETPOSTAVKA2}, we have the following
\begin{equation}\label{roUNula}
    \frac{e^{- \rho t}}{\rho e^{i \pi} \left( 1 + \frac{\lambda}{\rho ^{\rho e^{i \pi} A \left( \rho e^{i \pi} \right)} \left( e^{i \pi} \right)^{\rho e^{i \pi} A \left( \rho e^{i \pi} \right)}} \right)} \sim \frac{- 1}{\rho + \lambda e^{-i \pi \tilde{\alpha}} \rho^{1 - \tilde{\alpha}}}, \quad \rho \to 0.
\end{equation}
Using basic knowledge about convergence of integrals, if there are no other singularities on the negative part of the real axis, conditions \eqref{roUBeskonacno} and \eqref{roUNula} together guarantee that the integrals over $\Gamma_{h^{+}}$ and $\Gamma_{h^{-}}$ converge. If there were some singularities on the negative part of the real axis, we could modify our contour so that the branch cut is not done along the negative part of the real axis, but along some other ray passing through the origin of the complex plane.\\

Finally, summarizing all the calculations carried out, we arrive at the formula for the solution to equation \eqref{ScarpiRelaxation} in the general case
%\begin{equation}\label{konacnoResenje}
%\begin{aligned}
   %  &u(t) = \frac{u_{0}}{2 \pi i} \left( \mathfrak{R} + \int\limits_{0}^{\infty} \frac{\rho^{- \rho A \left( \rho e^{- i \pi} \right) - 1} e^{i \pi \left( \rho A \left( \rho e^{- i \pi} \right) + 1 \right)}}{\rho^{- \rho A \left( \rho e^{- i \pi} \right)} e^{i \pi \rho A \left( \rho e^{- i \pi} \right)} + \lambda} e^{- \rho t} d\rho \right. \\ &{} \hspace{13em}\left. - \int\limits_{0}^{\infty} \frac{\rho^{- \rho A \left( \rho e^{i \pi} \right) - 1} e^{- i \pi \left( \rho A \left( \rho e^{i \pi} \right) + 1 \right)}}{\rho^{- \rho A \left( \rho e^{i \pi} \right)} e^{- i \pi \rho A \left( \rho e^{i \pi} \right)} + \lambda} e^{- \rho t} d\rho \right).
%\end{aligned}\end{equation}
\begin{equation}
\begin{aligned}
     &u(t) = \frac{u_{0}}{2 \pi i} \left( \mathfrak{R} + \int\limits_{0}^{\infty} \frac{\rho^{- \rho A \left( \rho e^{- i \pi} \right) - 1} e^{i \pi \left( \rho A \left( \rho e^{- i \pi} \right) + 1 \right)}}{\rho^{- \rho A \left( \rho e^{- i \pi} \right)} e^{i \pi \rho A \left( \rho e^{- i \pi} \right)} + \lambda} e^{- \rho t} d\rho \right. \\ &{} \hspace{13em}\left. - \int\limits_{0}^{\infty} \frac{\rho^{- \rho A \left( \rho e^{i \pi} \right) - 1} e^{- i \pi \left( \rho A \left( \rho e^{i \pi} \right) + 1 \right)}}{\rho^{- \rho A \left( \rho e^{i \pi} \right)} e^{- i \pi \rho A \left( \rho e^{i \pi} \right)} + \lambda} e^{- \rho t} d\rho \right).
\end{aligned}\end{equation}
This ends our proof. $\square$
%Such a line must exist, consider the contrary, that for every segment with one endpoint at the origin and other endpoint in the left half plane, our function had a singularity, it would imply that the set of singularities has at least one accumulation point, and thus the apparatus of residual calculus would fall apart. {\sc Ne znam da li je ovaj komentar suvi\v san}.

%If we assume that the function $A(s)$ takes the same values for for $s_{1} = \rho e^{i \pi}$ and $s_{2} = \rho e^{- i \pi}$, for every positive real number $\rho$, we can then add the two integrals in \eqref{horizontalniIntegrali} and after some calculations we arrive at
%\begin{equation*}
 % \int\limits_{\Gamma_{h^{+}}} \frac{s^{s A(s) - 1}}{s^{s A(s)} + \lambda} e^{st} ds + \int\limits_{\Gamma_{h^{-}}} \frac{s^{s A(s) - 1}}{s^{s A(s)} + \lambda} e^{st} ds = \int\limits_{0}^{\infty} \frac{2 i \lambda \sin \left( \pi \rho A(- \rho) \right) e^{- \rho t}}{\rho^{- 2 \rho A(- \rho)} + \lambda^{2} + 2 \lambda \cos \left( \pi \rho A(- \rho) \right)} d\rho,
%\end{equation*}

\section{Special cases for transition function $\alpha(t)$}

In this section our aim is to apply general theory developed in the previous section to special cases of transition functions. As already said in Section 2, only functions satisfying certain properties will be among candidates to be considered as transition functions. Namely, we work with function $\alpha(t)$ whose LT exists and is explicitly known, and whose behaviour at the boundary points of its domain is controlled: $\lim\limits_{t\to0+}\alpha(t),\lim\limits_{t\to\infty}\alpha(t)\in(0,1)$. One of the main functions satisfying such conditions are certainly the exponential function with prescribed rate and the order transition function of Mittag-Leffler type (\cite{PERSPEKTIVA}, \cite{GarrappaNumerickoResenje}). It is also interesting to look at how our derived results, mainly \eqref{konacnoResenje} compare to well known results when the order of fractional derivative is constant. 

\subsection{Exponential transition}

First, assume that a transition function is given by
\begin{equation}\label{ekspTransition}
    \alpha(t) = \alpha_{2}+(\alpha_1-\alpha_2)e^{-ct}
\end{equation}
where $\alpha_{1}, \alpha_{2} \in (0, 1),\ c>0$. It is easy to see that a function of this type satisfies aforementioned conditions. This function describes a variable order transition from $\alpha_1$ to $\alpha_2$ obeying the exponential law with rate $-c$, and the final state is reached only asymptotically. However, we emphasize that such a choice for the transition function is no coincidence at all. Namely, physical nature of the transition function should be in compliance with physical nature of the solution we search for. Therefore, it is reasonable to assume that $\alpha(t)$ satisfies itself an evolution equation of the same type as \eqref{ScarpiRelaxation}. For the choice of $\alpha(t)$ as in \eqref{ekspTransition}, this is truly the case. Namely, the following system is satisfied: 
\begin{equation}
\begin{split}
\alpha'(t)=-c\alpha(t)+c\alpha_2,\\
\alpha(0)=\alpha_1.
\end{split}
\end{equation}
Preliminary numerical investigations with this choice of a transition function have been carried out by Garrappa et al. in \cite{GarrappaNumerickoResenje}. After analysis carried out in the previous section, however, we are able to provide an explicit form of the solution. Assume, for simplicity, that function \eqref{ekspTransition} is strictly increasing, that is $\alpha_1<\alpha_2$. We will prove the following:

\begin{Propozicija}\label{PROPOZICIJAZAEKSP}Let $u_0$ in \eqref{ScarpiRelaxation} be given and  assume that $\alpha(t)$ is an exponential transition function given by \eqref{ekspTransition}. Then, the unique solution of \eqref{ScarpiRelaxation} is given by
\begin{equation}\label{konacnoEksp}
    u(t) = \frac{u_{0}}{2\pi i} \left( \mathfrak{R} + \int\limits_{0}^{\infty} \frac{2 \lambda i \rho^{\frac{\alpha_{2} c - \alpha_{1} \rho}{c - \rho}-1} \sin \left( \pi \frac{\alpha_{2} c - \alpha_{1} \rho}{c - \rho} \right) e^{- \rho t}}{\rho^{2 \frac{\alpha_{2} c - \alpha_{1} \rho}{c - \rho}} + 2 \lambda \rho^{\frac{\alpha_{2} c - \alpha_{1} \rho}{c - \rho}} \cos \left( \pi \frac{\alpha_{2} c - \alpha_{1} \rho}{c - \rho} \right) + \lambda^{2}} d\rho \right),
\end{equation}
where $\mathfrak{R}$ denotes the sum of all residues of $\mathcal{L}[u](s)$ which are located inside $\Gamma$ multiplied with $2\pi i$, and $\Gamma$ is a contour sketched above.
\end{Propozicija}
\textbf{Proof: }One can clearly see that in that case the following holds:
\begin{equation*}
    A(s) = \frac{\alpha_{1} s + \alpha_{2} c}{s (s + c)}, \quad s^{s A(s)} = s^{\frac{\alpha_{1} s + \alpha_{2} c}{s + c}}, \quad s^{s A(s) - 1} = s^{\frac{(\alpha_{1} - 1) s + (\alpha_{2} - 1) c}{s + c}},\quad \mathit{Re}(s)>0.
\end{equation*}
We carry out our analysis with considerations above. Relation \eqref{TRAZENO} in this case becomes
\begin{equation}\label{DOBIJENO}\begin{aligned}
    \lim\limits_{\substack{R \to \infty \\ \varepsilon \to 0+}} \left( \int\limits_{\Gamma_{v}} \frac{s^{\frac{(\alpha_{1} - 1) s + (\alpha_{2} - 1) c}{s + c}}}{s^{\frac{\alpha_{1} s + \alpha_{2} c}{s + c} }+\lambda} e^{st} ds + \int\limits_{\Gamma_{h^+}} \frac{s^{\frac{(\alpha_{1} - 1) s + (\alpha_{2} - 1) c}{s + c}}}{s^{\frac{\alpha_{1} s + \alpha_{2} c}{s + c} }+\lambda} e^{st} ds + \int\limits_{\Gamma_{h^-}} \frac{s^{\frac{(\alpha_{1} - 1) s + (\alpha_{2} - 1) c}{s + c}}}{s^{\frac{\alpha_{1} s + \alpha_{2} c}{s + c} }+\lambda} e^{st} ds \right) = \mathfrak{R},
\end{aligned}\end{equation}
where $\mathfrak{R}$ denotes the sum of all residues inside the contour of integration multiplied by $2\pi i$. Singularities of function $\tilde{u}(s)e^{st}$ are $s_1=0,\,s_2=\infty,\,s_3=-c$, none of which lie in the contour of integration, and thus their residues do not count towards the solution. Remaining singularities are all solutions to the equation 
\begin{equation}\label{neznaniSingulariteti}
    s^{\frac{\alpha_{1} s + \alpha_{2} c}{s + c}} + \lambda = 0.
\end{equation}
Whether equation \eqref{neznaniSingulariteti} has any solutions in the complex plane, or what they look like, and most importantly whether they lie inside the contour of integration, remains an open problem. All that is known so far is that the equation \eqref{neznaniSingulariteti} has no solutions on the positive part of the real number line.
%Writing out equation \eqref{neznaniSingulariteti} in the polar form $s = \rho e^{i \theta}$, one obtains, after some calculations, that the following must hold for the solutions of \eqref{neznaniSingulariteti}:
%\begin{equation*}\begin{aligned}
    %&\frac{\left( \alpha_{1} \rho^{2} + \alpha_{2} c^{2} + \left( \alpha_{1} + \alpha_{2} \right) \rho c \cos \theta \right) \ln \rho - \left( \alpha_{1} - \alpha_{2} \right) \rho c \theta \sin \theta}{\rho^{2} + c^{2} + 2 \rho c \cos \theta} = \ln \lambda \\
    %&\frac{\left( \alpha_{1} \rho^{2} + \alpha_{2} c^{2} + \left( \alpha_{1} + \alpha_{2} \right) \rho c \cos \theta \right) \theta + \left( \alpha_{1} - \alpha_{2} \right) \rho c \ln \rho \sin \theta}{\rho^{2} + c^{2} + 2 \rho c \cos \theta} = \left( 2 k + 1 \right) \pi, \quad k \in \mathbb Z
%\end{aligned}\end{equation*}
%and from basic knowledge of complex analysis one knows that 
%\begin{equation}\label{SUMICA}
 %   \mathit{Res}[\tilde{u}(s)e^{st};s_1]+\mathit{Res}[\tilde{u}(s)e^{st};s_2]+\mathit{Res}[\tilde{u}(s)e^{st};s_3]+\mathit{Res}[\tilde{u}(s)e^{st};s_4]=0
%\end{equation}
%It is straightforward to calculate that $s_3$ and are points of removable singularities, thus their residues are equal to 0. After some lengthy calculation and with the help of the L'Hopital rule one gets that $s_1$ is a removable singularity as well, thus $\mathit{Res}[\tilde{u}(s)e^{st};s_1]=0$. Finally, $s_2$ appears to be a pole and by \eqref{SUMICA} we have that $\mathit{Res}[\tilde{u}(s)e^{st};s_2]=0$. 
In the light of \eqref{RESENJE} we have that
\begin{equation}\label{RESENJEKONKRETNO}
u(t)=\frac{u_{0}}{2\pi i}  \lim\limits_{\substack{R \to \infty \\ \varepsilon \to 0+}} \left( \mathfrak{R} -  \int\limits_{\Gamma_{h^+}} \frac{s^{\frac{(\alpha_{1} - 1) s + (\alpha_{2} - 1) c}{s + c}}}{s^{\frac{\alpha_{1} s + \alpha_{2} c}{s + c} }+\lambda} e^{st} ds - \int\limits_{\Gamma_{h^-}} \frac{s^{\frac{(\alpha_{1} - 1) s + (\alpha_{2} - 1) c}{s + c}}}{s^{\frac{\alpha_{1} s + \alpha_{2} c}{s + c} }+\lambda} e^{st} ds \right)   .
\end{equation}
Now, by definition of segments $\Gamma_{h^+}$ and $\Gamma_{h^-}$ one finds that
\begin{equation}\label{RESENJEVRLOKONKRETNO}
u(t)
=\frac{u_{0}}{2\pi i} \left( \mathfrak{R} - \int\limits_0^\infty\left( \frac{(\rho e^{i\pi})^{\frac{(\alpha_{1} - 1) \rho e^{i \pi} + (\alpha_{2} - 1) c}{\rho e^{i\pi} + c}}}{(\rho e^{i\pi})^{\frac{\alpha_{1} \rho e^{i \pi} + \alpha_{2} c}{\rho e^{i\pi} + c} }+\lambda}-\frac{(\rho e^{-i\pi})^{\frac{(\alpha_{1} - 1) \rho e^{-i \pi} + (\alpha_{2} - 1) c}{\rho e^{-i\pi} + c}}}{(\rho e^{-i\pi})^{\frac{\alpha_{1} \rho e^{-i \pi} + \alpha_{2} c}{\rho e^{-i\pi} + c} }+\lambda}\right)e^{-\rho t}d\rho \right)\end{equation}
With the analysis carried out in Section 2, we know that singularities $\rho = 0$ and $\rho = \infty$ are integrable. For given integrals to converge, we need to confirm that singularity $\rho = c$ is integrable when approached from the real axis. More specifically, we check the one-sided limits of the function under the integral sign as $\rho \to c$ along the real axis. First, using basic calculations, we rewrite the function in the following form:
\begin{equation}\label{drugiZapis}
    \frac{(\rho e^{i\pi})^{\frac{(\alpha_{1} - 1) \rho e^{i \pi} + (\alpha_{2} - 1) c}{\rho e^{i\pi} + c}}}{(\rho e^{i\pi})^{\frac{\alpha_{1} \rho e^{i \pi} + \alpha_{2} c}{\rho e^{i\pi} + c} }+\lambda} e^{- \rho t} = \frac{e^{- \rho t}}{- \rho \left( 1 + \lambda \rho^{\frac{\alpha_{1} \rho - \alpha_{2} c}{c - \rho}} e^{i \pi \frac{\alpha_{1} \rho - \alpha_{2} c}{c - \rho}} \right)}
\end{equation}
Now, recalling our assumption that $\alpha_{1} < \alpha_{2}$, when letting $\rho \to c^{-}$ we conclude that $\frac{\alpha_{1} \rho - \alpha_{2} c}{c - \rho} \to - \infty$ and that $\rho^{\frac{\alpha_{1} \rho - \alpha_{2} c}{c - \rho}} \to 0$. Also, one can easily see that $ e^{i \pi \frac{\alpha_{1} \rho - \alpha_{2} c}{c - \rho}}$ remains on the unit circle for all values of $\rho$, which, in turn, gives us $\lambda \rho^{\frac{\alpha_{1} \rho - \alpha_{2} c}{c - \rho}} e^{i \pi \frac{\alpha_{1} \rho - \alpha_{2} c}{c - \rho}} \to 0$ as $\rho \to c^{-}$, that is 
\begin{equation*}
    \lim_{\rho \to c^{-}} \frac{e^{- \rho t}}{- \rho \left( 1 + \lambda \rho^{\frac{\alpha_{1} \rho - \alpha_{2} c}{c - \rho}} e^{i \pi \frac{\alpha_{1} \rho - \alpha_{2} c}{c - \rho}} \right)} = \frac{e^{- c t}}{- c},
\end{equation*}
and this is a finite value. Similarly, taking limit $\rho \to c^{+}$ one deduces that $\frac{\alpha_{1} \rho - \alpha_{2} c}{c - \rho} \to \infty$ and consequently $\rho^{\frac{\alpha_{1} \rho - \alpha_{2} c}{c - \rho}} \to \infty$. Thus, the complex modulus of the denominator tends off to infinity as $\rho \to c^{+}$, which in turn, implies
\begin{equation*}
    \lim_{\rho \to c^{+}} \frac{e^{- \rho t}}{- \rho \left( 1 + \lambda \rho^{\frac{\alpha_{1} \rho - \alpha_{2} c}{c - \rho}} e^{i \pi \frac{\alpha_{1} \rho - \alpha_{2} c}{c - \rho}} \right)} = 0.
\end{equation*}
Given that one-sided limits of the function under the integral sign as $\rho \to c$ are finite, we conclude that the integrals in \eqref{RESENJEVRLOKONKRETNO} are convergent. Finally, combining expressions \eqref{drugiZapis} and \eqref{RESENJEVRLOKONKRETNO} one gets that the solution can be written more concisely as
\begin{equation}\label{konciznoResenje}
    u(t) = \frac{u_{0}}{2\pi i} \left( \mathfrak{R} - \int\limits_{0}^{\infty} \left( \frac{e^{- \rho t}}{\rho \left( 1 + \lambda \rho^{\frac{\alpha_{1} \rho - \alpha_{2} c}{c - \rho}} e^{- i \pi \frac{\alpha_{1} \rho - \alpha_{2} c}{c - \rho}} \right)} - \frac{e^{- \rho t}}{\rho \left( 1 + \lambda \rho^{\frac{\alpha_{1} \rho - \alpha_{2} c}{c - \rho}} e^{i \pi \frac{\alpha_{1} \rho - \alpha_{2} c}{c - \rho}} \right)} \right) d\rho \right),
\end{equation}
and after adding two fractions inside the integral and performing some calculations, using the definitions of sine and cosine in terms of complex exponentials, we arrive at
\begin{equation*}
    u(t) = \frac{u_{0}}{2\pi i} \left( \mathfrak{R} - \int\limits_{0}^{\infty} \frac{2 \lambda i \rho^{\frac{\alpha_{1} \rho - \alpha_{2} c}{c - \rho}-1} \sin \left( \pi \frac{\alpha_{1} \rho - \alpha_{2} c}{c - \rho} \right) e^{- \rho t}}{1 + 2 \lambda \rho^{\frac{\alpha_{1} \rho - \alpha_{2} c}{c - \rho}} \cos \left( \pi \frac{\alpha_{1} \rho - \alpha_{2} c}{c - \rho} \right) + \lambda^{2} \rho^{2 \frac{\alpha_{1} \rho - \alpha_{2} c}{c - \rho}}} d\rho \right).
\end{equation*}
After multiplying the denominator and numerator with $\rho^{2 \frac{\alpha_{2} c - \alpha_{1} \rho}{c - \rho}}$, and doing some basic algebra, we derive the final formula for the solution:
\begin{equation*}
     u(t) = \frac{u_{0}}{2\pi i} \left( \mathfrak{R} + \int\limits_{0}^{\infty} \frac{2 \lambda i \rho^{\frac{\alpha_{2} c - \alpha_{1} \rho}{c - \rho}-1} \sin \left( \pi \frac{\alpha_{2} c - \alpha_{1} \rho}{c - \rho} \right) e^{- \rho t}}{\rho^{2 \frac{\alpha_{2} c - \alpha_{1} \rho}{c - \rho}} + 2 \lambda \rho^{\frac{\alpha_{2} c - \alpha_{1} \rho}{c - \rho}} \cos \left( \pi \frac{\alpha_{2} c - \alpha_{1} \rho}{c - \rho} \right) + \lambda^{2}} d\rho \right).\ \square
\end{equation*}

\subsection{Order transition of Mittag-Leffler type}
Next, assume that a transition function is given by
\begin{equation}\label{MLTransition}
    \alpha(t) = \alpha_{2}+(\alpha_1-\alpha_2)E_\beta(-ct^\beta)
\end{equation}
where $\alpha_{1}, \alpha_{2} \in (0, 1)$ and $$E_\beta(z)=\sum\limits_{k=0}^\infty\frac{z^k}{\Gamma(\alpha k+\beta)}$$ is the one parameter ML function. It is easy to see that a function of this type also satisfies conditions stated in subsection 2.1. For the choice of $\alpha(t)$ as in \eqref{MLTransition}, the following system is satisfied: 
\begin{equation}
\begin{split}
\prescript{C}{}{D}^\beta u(t)=-c\alpha(t)+c\alpha_2,\\
\alpha(0)=\alpha_1.
\end{split}
\end{equation}
Assume, as in the previous subsection, that function \eqref{MLTransition} is strictly increasing, that is $\alpha_1<\alpha_2$. One can clearly see that in that case the following holds:
\begin{equation*}
    A(s) = \frac{\alpha_{1} s^\beta + \alpha_{2} c}{s (s^\beta + c)}, \quad s^{s A(s)} = s^{\frac{\alpha_{1} s^\beta + \alpha_{2} c}{s^\beta + c}}, \quad s^{s A(s) - 1} = s^{\frac{(\alpha_{1} - 1) s^\beta + (\alpha_{2} - 1) c}{s^\beta + c}},\quad \vert s\vert>\vert c\vert^{1/\beta}.
\end{equation*}
Using the similar calculus as in the previous subsection, one can obtain the following theorem.

\begin{Propozicija}Let $u_0$ in \eqref{ScarpiRelaxation} be given and  assume that $\alpha(t)$ is a Mittag-Leffler transition function given by \eqref{MLTransition}. Then, the unique solution of \eqref{ScarpiRelaxation} is given by
\begin{equation}\label{RESENJEVRLOKONKRETNOML}
u(t)
=\frac{u_{0}}{2\pi i} \left( \mathfrak{R} - \int\limits_0^\infty\left( \frac{(\rho e^{i\pi})^{\frac{(\alpha_{1} - 1) (\rho e^{i \pi})^\beta + (\alpha_{2} - 1) c}{(\rho e^{i\pi})^\beta + c}}}{(\rho e^{i\pi})^{\frac{\alpha_{1} (\rho e^{i \pi})^\beta + \alpha_{2} c}{(\rho e^{i\pi})^\beta + c} }+\lambda}-\frac{(\rho e^{-i\pi})^{\frac{(\alpha_{1} - 1) (\rho e^{-i \pi})^\beta + (\alpha_{2} - 1) c}{(\rho e^{-i\pi})^\beta + c}}}{(\rho e^{-i\pi})^{\frac{\alpha_{1} (\rho e^{-i \pi})^\beta + \alpha_{2} c}{(\rho e^{-i\pi})^\beta + c} }+\lambda}\right)e^{-\rho t}d\rho \right)\end{equation}
where $\mathfrak{R}$ denotes the sum of all residues of $\mathcal{L}[u](s)$ which are located inside $\Gamma$ multiplied with $2\pi i$, and $\Gamma$ is a contour sketched above.
\end{Propozicija}
\textbf{Proof: }As already said, $s^{s A(s)} = s^{\frac{\alpha_{1} s^\beta + \alpha_{2} c}{s^\beta + c}}$, and $s^{s A(s) - 1} = s^{\frac{(\alpha_{1} - 1) s^\beta + (\alpha_{2} - 1) c}{s^\beta + c}}$ for $\vert s\vert>\vert c\vert^{1/\beta}.$ Having this in mind, we derive expression \eqref{RESENJEVRLOKONKRETNOML} exactly in the same way as we derived expression \eqref{RESENJEVRLOKONKRETNO} in the proof of Proposition \eqref{PROPOZICIJAZAEKSP}.

\subsection{Fractional derivative of constant order}

Now we turn to the case when we have a constant-order fractional derivative, namely when 
\begin{equation*}
    \alpha(t) = \alpha, \; \alpha \in (0, 1).
\end{equation*}
It is now clear that we have 
\begin{equation*}
    A(s) = \frac{\alpha}{s}, \quad s^{s A(s)} = s^{\alpha}, \quad s^{s A(s) - 1} = s^{\alpha - 1}.
\end{equation*}
and thus the solution to problem \eqref{ScarpiRelaxation} reduces to finding the following inverse:
\begin{equation}\label{konstLambda}
    u(t) = u_{0} \mathcal{L}^{-1} \left[ \frac{s^{\alpha - 1}}{s^{\alpha} + \lambda} \right](t).
\end{equation}
In \cite{IntMittagLeffler} F. Mainardi defines $e_{\alpha} (t) : = E_{\alpha} (- t^{\alpha})$ such that 
\begin{equation*}
    e_{\alpha}(t) = \mathcal{L}^{-1} \left[ \frac{s^{\alpha - 1}}{s^{\alpha} + 1} \right](t)
\end{equation*}
holds for $\alpha > 0$, while also giving the following integral representation
\begin{equation*}
    e_{\alpha}(t) = \frac{1}{\pi} \int\limits_{0}^{\infty} \frac{r^{\alpha - 1} \sin (\alpha \pi)}{r^{2 \alpha} + 2 r^{\alpha} \cos(\alpha \pi) + 1} e^{- r t} dr
\end{equation*}
which holds for $0 < \alpha < 1$. Using simple calculations and basic properties of LT, we arrive at
\begin{equation*}
    \mathcal{L}^{-1} \left[ \frac{s^{\alpha - 1}}{s^{\alpha} + \lambda} \right](t) = \mathcal{L}^{-1} \left[\frac{1}{\lambda^{\frac{1}{\alpha}}} \frac{\left( \frac{s}{\lambda^{\frac{1}{\alpha}}} \right)^{\alpha - 1}}{\left( \frac{s}{\lambda^{\frac{1}{\alpha}}} \right)^{\alpha} + 1} \right](t) = e_{\alpha} \left( \lambda^{\frac{1}{\alpha}} t \right).
\end{equation*}
Therefore, \eqref{konstLambda} simply becomes
\begin{equation*}
    u(t) = u_{0} e_{\alpha} \left( \lambda^{\frac{1}{\alpha}} t \right) = u_{0} \frac{1}{\pi} \int\limits_{0}^{\infty} \frac{r^{\alpha - 1} \sin (\alpha \pi)}{r^{2 \alpha} + 2 r^{\alpha} \cos(\alpha \pi) + 1} e^{- r \lambda^{\frac{1}{\alpha}}t} dr
\end{equation*}
and with a standard change of variables $\rho = r \lambda^{\frac{1}{\alpha}}$ in the integral above, we arrive at the solution
\begin{equation}\label{konstInt}
    u(t) = u_{0} \frac{1}{\pi} \int\limits_{0}^{\infty} \frac{\lambda \rho^{\alpha - 1} \sin (\alpha \pi)}{\rho^{2 \alpha} + 2 \lambda \rho^{\alpha} \cos (\alpha \pi) + \lambda^{2}} e^{- \rho t} d\rho. 
\end{equation}
Now, we turn to formula \eqref{konacnoResenje}. Knowing that $\displaystyle A(s) = \frac{\alpha}{s}$ we first have
\begin{equation*}
    u(t) = \frac{u_{0}}{2 \pi i} \left( \mathfrak{R} + \int\limits_{0}^{\infty} \left( \frac{\rho^{-\rho \frac{\alpha}{\rho e^{- i \pi}}- 1} e^{i \pi \left( \rho \frac{\alpha}{\rho e^{- i \pi}}+ 1 \right)}}{\rho^{-\rho \frac{\alpha}{\rho e^{- i \pi}}}e^{i \pi\rho \frac{\alpha}{\rho e^{- i \pi}}} + \lambda} - \frac{\rho^{-\rho \frac{\alpha}{\rho e^{ i \pi}}- 1} e^{-i \pi \left( \rho \frac{\alpha}{\rho e^{ i \pi}}+ 1 \right)}}{\rho^{-\rho \frac{\alpha}{\rho e^{ i \pi}}}e^{-i \pi\rho \frac{\alpha}{\rho e^{ i \pi}}} + \lambda}  \right) e^{- \rho t} d\rho \right).
\end{equation*}
After performing some calculations and using the definition of complex exponentials we arrive at
\begin{equation*}
    u(t) = \frac{u_{0}}{2 \pi i} \left( \mathfrak{R} + \int\limits_{0}^{\infty} \frac{2i \lambda \rho^{\alpha -1} \sin(\alpha \pi)}{\rho^{2 \alpha} + 2 \lambda \rho^{\alpha} \cos(\alpha \pi) + \lambda^{2}} e^{- \rho t} d\rho \right)
\end{equation*}
It is known that the singularities of function $\displaystyle \frac{s^{\alpha - 1}}{s^{\alpha} + \lambda} e^{st}$ are exactly the solutions to $s^{\alpha} + \lambda = 0$, that is $s = \lambda^{\frac{1}{\alpha}} e^{\frac{\left( 2k + 1 \right) \pi i}{\alpha}}$ where $k \in \mathbb Z$, but none of these singularities have argument in $(- \pi, \pi)$ because $0 < \alpha < 1$, and thus they do not contribute to $\mathfrak{R}$. Therefore, we finally arrive at 
\begin{equation*}
    u(t) = \frac{u_{0}}{\pi} \int\limits_{0}^{\infty} \frac{\lambda \rho^{\alpha -1} \sin(\alpha \pi)}{\rho^{2 \alpha} + 2 \lambda \rho^{\alpha} \cos(\alpha \pi) + \lambda^{2}} e^{- \rho t} d\rho
\end{equation*}
which is exactly \eqref{konstInt}. Similar calculations can be performed by setting $c = 0$ in \eqref{konacnoEksp} to get \eqref{konstInt} with $\alpha_{1}$ in place of $\alpha$, or by letting $c \to \infty$ with $\alpha_{2}$ in place of $\alpha$. We conclude that the results we derived generalize well known results for the constant-order fractional derivative.

\section{Numerical experiments}
In this section we employ certain numerical methods to show reliability of previously stated results. This will be done by plotting the solution $u(t)$ of the relaxation equation \eqref{ScarpiRelaxation} for various choices of transition function $\alpha(t)$ and comparing it with the solutions $u_1(t)$ and $u_2(t)$ of corresponding fractional relaxation equations of orders $\alpha_1$ and $\alpha_2$. One can notice that we cannot plot the graph of $u(t)$ with the greatest precision, since the expression for $u(t)$ (see \eqref{konacnoResenje}) includes an unknown variable $\mathfrak{R}$. However, if we write that solution in the form $u(t)=\frac{u_0}{2\pi\imath}\mathfrak{R}+v(t)$, where real function $v(t)$ is calculated via indefinite integrals in \eqref{konacnoResenje}, then it is enough to conclude that an appropriate translate of the graph of $v(t)$ approximates graphs of $u_1(t)$ and $u_2(t)$. We choose to do this for the choice of an exponential transition function, since the choice of Mittag-Leffler transition offers no numerical novelty but more complicated computing solely.

In the sequel we show a few sketches where graphs of $u_1(t)$ and $u_2(t)$ are obtained via definition of Mittag-Leffler function, and the graph of $u(t)$ is found using the Gauss-Kronrad quadrature formula. We provide a fragment of an appropriate MATLAB code at the end of this section.

\begin{figure}[H]\label{SLIKAZANUMERIKU}
\centering
\includegraphics[width=0.63\linewidth,trim=0 540 250 0, clip]{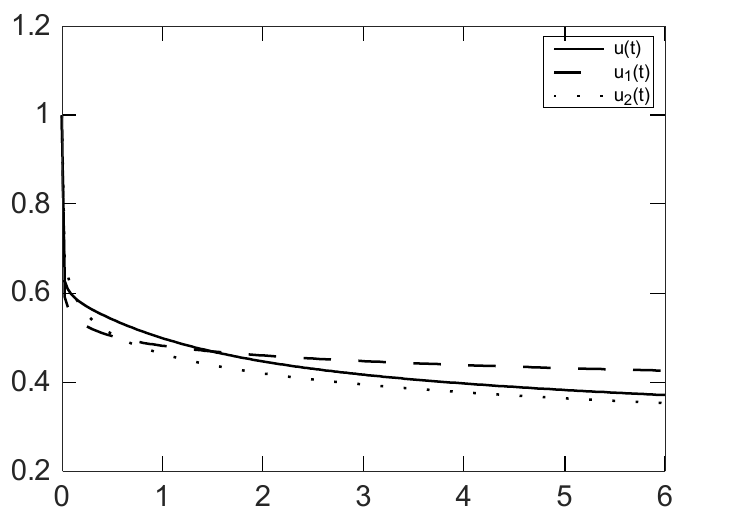}
\caption{Plot of the translated solution $u(t)$ of the relaxation equation \eqref{ScarpiRelaxation} with $u_0=\lambda=1$ for exponential transition $\alpha(t)=\alpha_2+(\alpha_1-\alpha_2)e^{-ct}$ with $\alpha_1=\frac{1}{8},\alpha_2=\frac{1}{4}$ and $c=2$, and of the solutions $u_1(t)$ and $u_2(t)$ of constant-order relaxation equations of orders $\alpha_1$ and $\alpha_2$.}
\end{figure}
Observe that in Figure 2 solution $u(t)$ approaches $u_1(t)$ as $t\to0$ and $u(t)$ approaches $u_2(t)$ as $t\to\infty$. This phenomenon is in compliance with \cite[Figure 2, p.9]{GarrappaNumerickoResenje}, and is a certain confirmation of \cite[Proposition 5]{GarrappaNumerickoResenje}. Observe also that $u(0)=u_1(0)=u_2(0)=1$ for our choice of parameters, and this is also visible in Figure 2. In the sequel we give two more figures (Figure 3 and Figure 4) obtained for various choices of $\alpha_1,\alpha_2$ and $c$.

\begin{figure}[H]\label{SLIKAZANUMERIKU2}
\centering
\includegraphics[width=0.63\linewidth,trim=0 540 250 0, clip]{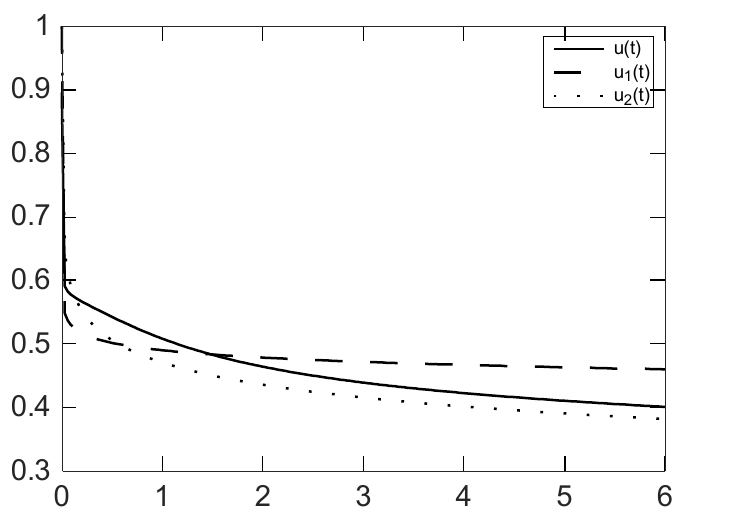}
\caption{Plot of the translated solution $u(t)$ of the relaxation equation \eqref{ScarpiRelaxation} with $u_0=\lambda=1$ for exponential transition $\alpha(t)=\alpha_2+(\alpha_1-\alpha_2)e^{-ct}$ with $\alpha_1=\frac{1}{15},\alpha_2=\frac{1}{5}$ and $c=2$, and of the solutions $u_1(t)$ and $u_2(t)$ of constant-order relaxation equations of orders $\alpha_1$ and $\alpha_2$.}
\end{figure}

\begin{figure}[H]\label{SLIKAZANUMERIKU3}
\centering
\includegraphics[width=0.63\linewidth,trim=0 545 250 0, clip]{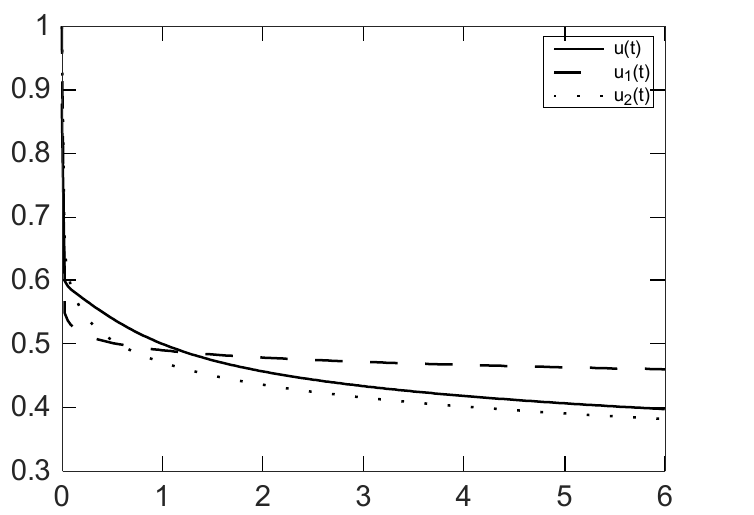}
\caption{Plot of the translated solution $u(t)$ of the relaxation equation \eqref{ScarpiRelaxation} with $u_0=\lambda=1$ for exponential transition $\alpha(t)=\alpha_2+(\alpha_1-\alpha_2)e^{-ct}$ with $\alpha_1=\frac{1}{15},\alpha_2=\frac{1}{5}$ and $c=3$, and of the solutions $u_1(t)$ and $u_2(t)$ of constant-order relaxation equations of orders $\alpha_1$ and $\alpha_2$.}
\end{figure}

Notice that the last two figures also confirm \cite[Proposition 5]{GarrappaNumerickoResenje} and are compatible with \cite[Figure 2, p.9]{GarrappaNumerickoResenje}.

Finally, we provide a fragment of a Matlab code that illustrates Gauss-Kronrad quadrature rule. In the following code Matlab function $integrand\_fun(x,t)$ is determined by the expression under integral sign in \eqref{konacnoEksp}.

\begin{center}
\doublebox{
\begin{minipage}{0.8\linewidth}
\begin{verbatim}
    % High-accuracy integration options
    absTol = 1e-12;
    relTol = 1e-12;

    % Define t grid
    t_values = linspace(0, 6, 200);

    % Compute u(t) for each t
    for k = 1:length(t_values)
        t = t_values(k);

        integrand = @(x) integrand_fun(x, t);


    [u_values(k),err]= quadgk(integrand, 0, Inf, ...
            'RelTol', relTol, ...
            'AbsTol', absTol, ...
            'ArrayValued', true);
    end
\end{verbatim}
\end{minipage}
}
\end{center}

\section*{Declarations}

{\bf Conflict of interests:} Authors declare that they have no conflicts of interest.

\section*{Data availability statement}

Authors declare that no datasets have been generated during this study.

\section*{Acknowledgements}

Authors are grateful to professors Ljubica Oparnica and Du\v{s}an Zorica for introduction to this topic and for useful discussions regarding research areas in this field.\\
We also express gratitude to anonymous referees whose suggestions improved the appearance of this manuscript.

\end{document}